# GLOBAL FLOWS FOR STOCHASTIC DIFFERENTIAL EQUATIONS WITHOUT GLOBAL LIPSCHITZ CONDITIONS

BY SHIZAN FANG, PETER IMKELLER AND TUSHENG ZHANG

*Université de Bourgogne, Humboldt-Universität zu Berlin and University of Manchester*

We consider stochastic differential equations driven by Wiener processes. The vector fields are supposed to satisfy only local Lipschitz conditions. The Lipschitz constants of the drift vector field, valid on balls of radius $R$, are supposed to grow not faster than $\log R$, while those of the diffusion vector fields are supposed to grow not faster than $\sqrt{\log R}$. We regularize the stochastic differential equations by associating with them approximating ordinary differential equations obtained by discretization of the increments of the Wiener process on small intervals. By showing that the flow associated with a regularized equation converges uniformly to the solution of the stochastic differential equation, we simultaneously establish the existence of a global flow for the stochastic equation under local Lipschitz conditions.

**Introduction.** Let $A_0, A_1, \ldots, A_N$ be $N+1$ vector fields on the Euclidean space $\mathbf{R}^d$ and $(w_t)_{t \geq 0}$ be an $\mathbf{R}^N$-valued standard Brownian motion. Consider the following Stratonovich stochastic differential equation:

$$(0.1) \qquad dx_t = \sum_{i=1}^N A_i(x_t) \circ dw_t^i + A_0(x_t)\,dt, \qquad x_0 = x,$$

where $w_t^i$ denotes the $i$th component of $w_t$. If the coefficients are sufficiently smooth, for example, if $A_1, \ldots, A_N$ are $\mathcal{C}^2$ and $A_0$ is $\mathcal{C}^1$, then the stochastic differential equation (0.1) has a unique solution $(x_t)$. In this case, $(x_t)$ solves also the following Itô stochastic differential equation:

$$(0.2) \qquad dx_t = \sum_{i=1}^{N} A_i(x_t)\,dw_t^i + \tilde{A}_0(x_t)\,dt, \qquad x_0 = x,$$

where

$$(0.3) \qquad \tilde{A}_0 = A_0 + \frac{1}{2}\sum_{i=1}^{N}\sum_{j=1}^{d}\frac{\partial A_i}{\partial x_j}A_i^j.$$

Under global Lipschitz conditions on the coefficients $\tilde{A}_0, A_1, \ldots, A_N$, Kunita [7] proved that the Itô stochastic differential equation (0.2) defines a global flow of homeomorphisms. On the other hand, subject to the assumptions that the coefficients $A_1, \ldots, A_N$ are $\mathcal{C}^2$, bounded and with bounded derivatives of first and second order, and that $A_0$ is $\mathcal{C}^1$, also bounded with bounded derivative, Moulinier [10] proved that almost surely the solutions $(x_t^n)$ of the regularized ordinary differential equations

$$(0.4) \qquad dx_t^n = \sum_{i=1}^{N} A_i(x_t^n)\dot{w}_t^{n,i}\,dt + A_0(x_t^n)\,dt, \qquad x_0^n = x,$$

where

$$(0.5) \quad \dot{w}_t^n = 2^n(w_{(k+1)2^{-n}} - w_{k2^{-n}}) \qquad \text{for } t \in [k2^{-n},(k+1)2^{-n}[,\ k\geq 0,$$

converge to $(x_t)$, uniformly with respect to $(t,x)$ in each compact subset of $\mathbf{R}_+ \times \mathbf{R}^d$. This provides another approach to the existence of global flows. For related work, we refer to [1, 2, 4, 9] and [11].

The main aim of this paper is to remove the global Lipschitz conditions from the hypotheses needed to arrive at these conclusions. Based on moment estimates for the one-point and two-point motions with explicit dependence on the Lipschitz constants, we still obtain the smooth approximation to the solution of (0.1). Consequently, we will prove the following result:

THEOREM A. *Let $A_1, \ldots, A_N$ be in the class $\mathcal{C}^2$ and $A_0$ in the class $\mathcal{C}^1$. Suppose that* (i) *the growth of the coefficients $A_1, \ldots, A_N$ and their first and second order derivatives is dominated by $\sqrt{\log|x|}$ and* (ii) *the growth of $A_0$ and its first order derivatives is controlled by $\log|x|$, as $|x| \to \infty$. Then the Stratonovich stochastic differential equation* (0.1) *defines a global flow of homeomorphisms, that is, for each $t > 0$, the map $x \to x_t(x)$ is a homeomorphism of $\mathbf{R}^d$ almost surely.*

The moment estimates we propose in order to derive this theorem are in the spirit of Imkeller and Scheutzow [5] and Imkeller [6]. Their starting point is a specification of the constant $c_p(p \geq 1)$ in an inequality of the type

$$\mathbf{E}\bigg(\sup_{0\leq t\leq 1}|x_t(x) - x_t(y)|^p\bigg) \leq c_p|x-y|^p, \qquad x,y \in \mathbf{R}^d.$$

In [5], under global Lipschitz conditions on the vector fields, it is seen to be essentially given by $c_p = \exp(cp^2)$, with a universal $c$. Here, we shall work with similar ideas. The essential novelty is the following observation. Assume that only local Lipschitz conditions are given, which on large balls of radius $m$ centered at the origin are given by $L_m$. For each $m$, replace the original vector fields by vector fields with global Lipschitz conditions and Lipschitz constant essentially equal to $L_m$. Then the global two-point motions $(x_t(x), x_t(y))$ are related to the two-point motions $(x_t^m(x), x_t^m(y))$ associated with the modified vector fields through the key equality

$$|x_t(x) - x_t(y)|^p = \sum_{m=1}^{\infty} |x_t^m(x) - x_t^m(y)|^p \mathbf{1}_{\{m-1 \leq Y_1(x) \vee Y_1(y) < m\}},$$

where $t \in [0,1]$ and $Y_1(x) = \sup_{0 \leq t \leq 1} |x_t(x)|, x \in \mathbf{R}^d$. So, the two-point motions of the global flow will be controlled by the two-point motions of the modified flows and the growth behavior of the one-point motions of the global flow. This idea is exploited in Section 1 below (Theorems 1.7 and 1.8). Section 2 is devoted to giving moment estimates of the same type for regularized ordinary differential equations obtained by discretizing the increments of the Wiener process on dyadic time intervals. Again, this is done for one- and two-point motions separately. But the discretization procedure will produce a bad term $e^{\alpha_n}$ (see Theorem 2.6), where $\alpha_n$ is an exponential function of Lipschitz constants. In order to get the desired result, in Section 3, we truncate the vector fields and simultaneously discretize the Wiener process. In this way, the bad term $e^{\alpha_n}$ can be handled. We show, using our moment inequality techniques, that the flows of the regularized ordinary differential equations converge to the flows of the original stochastic differential equation, under local Lipschitz conditions (Theorem 3.4). The Lipschitz constants on balls of radius $R$ centered at zero are of order $\log R$ for the drift vector field and $\sqrt{\log R}$ for the diffusion vector fields. These conditions constitute hypothesis (H). Hence, Theorem 3.4 implies Theorem A in the usual way (see, e.g., [4]).

We should mention that the existence of global flows of homeomorphisms for *one-dimensional* stochastic differential equations was established by Yamada and Ogura [12], under local Lipschitz and linear growth conditions on the coefficients. For the multidimensional case, the situation is quite different; in fact, if we denote by $\tau_x$ the lifetime of the solution $(x_t(x))$ to the stochastic differential equation (0.1), the linear growth (or even boundedness) of coefficients is not sufficient to ensure that

(0.6) $$P(\tau_x = +\infty, \text{ for all } x \in \mathbf{R}^d) = 1.$$

In the case where the diffusion coefficients are in $\mathcal{C}^{2+\delta}$ and the drift is $\mathcal{C}^{1+\delta}$ with $\delta > 0$, using local flows of derivatives of solutions, Li [8] proved (0.6) for the stochastic differential equation (0.1), as well as for its dual equation (see [7] for this notion), under the same growth condition on the local Lipschitz constants as ours in Theorem A; therefore, by Theorem 6.1 or Theorem 7.3 in [7], she obtains a global flow of diffeomorphisms. Note that even for Itô stochastic differential equations, smoothness of coefficients with $\delta > 0$ was needed to apply Theorem 6.1 of [7]. For a study of stochastic differential equations under non(local) Lipschitz conditions, we refer to [3].

**1. Moment estimates for one- and two-point motions.** Let $(x_t(x))$ be the solution of the Itô stochastic differential equation (0.2). The growth of the moments of $(x_t(x))$ in the spatial parameter will crucially depend on the growth behavior of the diffusion coefficients $A_1, \ldots, A_N$. In order to capture well the growth of the local Lipschitz constants for estimating moments of the two-point motions $\mathbf{E}(|x_t(x) - x_t(y)|^p)$, we shall distinguish between the following hypotheses:

(H1) there are constants $C_1$ and $C_2 > 0$ such that

$$\sum_{i=1}^{N} |A_i(x)|^2 \leq C_1^2, \qquad |\tilde{A}_0(x)| \leq C_2(1 + |x|);$$

(H2) there are constants $C_3$ and $C_4 > 0$ such that

$$\sum_{i=1}^{N} |A_i(x)|^2 \leq C_3^2(1 + |x|^2), \qquad |\tilde{A}_0(x)| \leq C_4(1 + |x|).$$

Let us remark at this point that our setting could be extended to the case of infinitely many vector fields and, correspondingly, an infinite-dimensional Wiener process, by noting that (H1) and (H2) only concern Euclidean norms and could be stated for Hilbert–Schmidt norms instead. In what follows, universal positive constants appearing in the inequalities are denoted by $C$ and allowed to change from instance to instance.

1.1. *Precise $L^p$-estimates for the one-point motion.* Define $Y_t(x) = \sup_{0 \leq s \leq t} |x_s(x)|$. We shall first give the explicit estimate of $\|Y_1(x)\|_p$ as a function of $p$.

PROPOSITION 1.1. *Under condition* (H1), *we have for any $p > 1$,*

(1.1) $$\|Y_1(x)\|_p \leq (1 + CC_1\sqrt{p})e^{C_2}(1 + |x|).$$

PROOF. Fix $x \in \mathbf{R}^d$. For $0 \leq t \leq 1$, put $\varphi(t) = \|Y_t(x)\|_p$ and $M_t = \sum_{i=1}^{N} \int_0^t A_i(x_s(x)) \, dw_s^i$. By the inequality of Burkholder, Davis and Gundy (see [5]), for any $0 \leq T \leq 1$,

$$\mathbf{E}\left(\sup_{0 \leq t \leq T} |M_t|^p\right) \leq C\sqrt{p^p}\mathbf{E}\left[\left(\int_0^T \sum_{i=1}^{N} |A_i(x_s(x))|^2 \, ds\right)^{p/2}\right] \leq CC_1^p\sqrt{p^p},$$

or
$$\left\| \sup_{0 \le t \le T} |M_t|^p \right\|_p \le CC_1 \sqrt{p}.$$

Using (0.2), we get the inequality
$$\varphi(T) \le |x| + CC_1\sqrt{p} + C_2 \int_0^T (1 + \varphi(s))\, ds.$$

Dividing both sides by the term $1 + |x|$ and applying Gronwall's lemma to the function $\varphi(T) + 1/(1 + |x|)$, we get $\frac{\varphi(1)}{1+|x|} \le (1 + CC_1\sqrt{p})e^{C_2}$ and estimate (1.1) follows. $\square$

The preceding moment inequality implies the following exponential inequality:

COROLLARY 1.2. *Suppose that* (H1) *holds. For $R > 0$, there exists $\delta_0 = \delta_0(C_1, C_2, R) > 0$ such that*

(1.2) $$\sup_{|x| \le R} \mathbf{E}(e^{\delta_0 Y_1^2(x)}) < +\infty.$$

PROOF. By (1.1), there is a constant $\beta$ such that $\|Y_1(x)\|_p \le \beta\sqrt{p}(1 + |x|)$. Let $\delta > 0$. We have

$$\mathbf{E}(e^{\delta Y_1^2(x)}) = 1 + \sum_{p=1}^{+\infty} \frac{\delta^p \mathbf{E}(Y_1^{2p}(x))}{p!} \le 1 + \sum_{p=1}^{+\infty} \frac{\delta^p \beta^{2p}(2p)^p(1+|x|)^{2p}}{p!}.$$

By Stirling's formula, $\frac{p^p}{p!} \sim \frac{e^p}{\sqrt{2\pi p}}$ as $p \to +\infty$. Hence, if $|x| \le R$, then the above expression is dominated by

$$C\left(1 + \sum_{p=1}^{+\infty} (2\delta\beta^2 e(1+R)^2)^p\right) = \frac{C}{1 - 2\delta\beta^2 e(1+R)^2},$$

which is finite if $\delta < 1/(2\beta^2 e(1+R)^2)$. So we get (1.2). $\square$

In the following proposition, we shall investigate estimates under (H2):

PROPOSITION 1.3. *Under condition* (H2), *there are constants $\beta_1$ and $\beta_2 > 0$ such that for all $p > 1$ and $x \in \mathbf{R}^d$,*

(1.3) $$\|Y_1(x)\|_p \le \beta_1 e^{\beta_2 p}(1 + |x|).$$

PROOF. Let $M$ and $\varphi$ be defined as in the proof of Proposition 1.1. Under (H2), we have

$$\left\|\sup_{0\leq t\leq T}|M_t|\right\|_p \leq CC_1\sqrt{p}\left[\int_0^T(1+\varphi^2(s))\,ds\right]^{1/2}.$$

Therefore, in this case, the inequality

$$(1.4) \quad \varphi(T) \leq |x| + CC_1\sqrt{p}\left[\int_0^T(1+\varphi^2(s))\,ds\right]^{1/2} + C_2\int_0^T(1+\varphi^2(s))\,ds$$

follows. To apply Gronwall's lemma, we must square the two sides of (1.4), which results in

$$\varphi^2(T) \leq 3\bigg(|x|^2 + (C^2C_1^2 p + 2C_2^2)\int_0^T(1+\varphi^2(s))\,ds\bigg).$$

It follows that

$$\frac{\varphi^2(T)+1}{(1+|x|)^2} \leq 3\exp\{3(C^2C_1^2 p + 2C_2^2)T\},$$

from which we deduce (1.3). $\square$

In the same spirit, we can treat the time variation of the one-point motion moments.

COROLLARY 1.4. *Under hypothesis* (H1) *or* (H2), *for any $p>1$, there exists a constant $C_p > 0$ (which depends on $C_1$ and $C_2$, or on $C_3$ and $C_4$, resp.) such that for $x \in \mathbf{R}^d, s, t \geq 0$,*

$$(1.5) \qquad \mathbf{E}(|x_t(x) - x_s(x)|^{2p}) \leq C_p |t-s|^p(1+|x|)^{2p}.$$

PROOF. We have, for $s < t$, $x \in \mathbf{R}^d$,

$$x_t(x) - x_s(x) = \sum_{i=1}^N \int_s^t A_i(x_u(x))\,dw_u^i + \int_s^t \tilde{A}_0(x_u(x))\,du.$$

Hence, there exists a constant $\beta_p > 0$ such that

$$\mathbf{E}(|x_t(x) - x_s(x)|^{2p}) \leq \beta_p\bigg\{\mathbf{E}\bigg[\bigg(\int_s^t \sum_{i=1}^N |A_i(x_u(x))|^2\,du\bigg)^p\bigg]$$
$$+ \mathbf{E}\bigg[\bigg(\int_s^t |\tilde{A}_0(x_u(x))|\,du\bigg)^{2p}\bigg]\bigg\}.$$

So, we see that for some sufficiently large constant $C_p > 0$, the right-hand side of the above inequality is dominated by

$$C_p(t-s)^p(1+\mathbf{E}(Y_1(x)^{2p})).$$

We now obtain (1.5) for an eventually different $C_p$ by using (1.1) or (1.3). $\square$

1.2. *Precise $L^p$-estimates for the two-point motion under global Lipschitz conditions.* Here, we shall work under the following global Lipschitz condition:

(L) there exist constants $L_1$ and $L_2 > 0$ such that

$$\sum_{i=1}^{N} |A_i(x) - A_i(y)|^2 \leq L_1^2 |x-y|^2,$$

$$|\tilde{A}_0(x) - \tilde{A}_0(y)| \leq L_2 |x-y|, \qquad x,y \in \mathbf{R}^d.$$

Set $Y_T(x,y) = \sup_{0 \leq t \leq T} |x_t(x) - x_t(y)|$. We shall give the explicit dependence on $L_1$ and $L_2$ for $L^p$-estimates of $Y_1(x,y)$.

PROPOSITION 1.5. *Under hypothesis* (L), *we have for any $p > 1$, all $x, y \in \mathbf{R}^d$,*

(1.6) $$\mathbf{E}(Y_1(x,y)^p) \leq 2^p |x-y|^p e^{C L_1^2 p^2 + L_2^2 p}.$$

PROOF. Let $\varphi(T) = \|Y_T(x,y)\|_p$. As in the estimates above, we have

(1.7) $$\varphi(T) \leq |x-y| + C L_1 \sqrt{p} \left[\int_0^T \varphi^2(s)\,ds\right]^{1/2} + L_2 \int_0^T \varphi(s)\,ds.$$

Squaring the two sides of (1.7) results in

$$\varphi^2(T) \leq 2\left(2|x-y|^2 + (2C^2 L_1^2 p + L_2^2) \int_0^T \varphi^2(s)\,ds\right), \qquad T \leq 1.$$

It follows that for an eventually different constant $C > 0$,

$$\varphi(1) \leq 2|x-y| e^{CL_1^2 p + L_2^2},$$

from which we obtain (1.6). □

REMARK. In squaring the two sides of (1.7), the control on the Lipschitz constant $L_2$ was lost. In order to recapture it, we shall now consider only $\mathbf{E}(|x_t(x) - x_t(y)|^p)$.

PROPOSITION 1.6. *Assume* (L). *Then for any $p \geq 2$, all $x, y \in \mathbf{R}^d$ and $t \in [0,1]$,*

(1.8) $$\mathbf{E}(|x_t(x) - x_t(y)|^{2p}) \leq |x-y|^{2p} e^{2p^2 L_1^2 + 2pL_2}.$$

PROOF. Let $\xi_t = |x_t(x) - x_t(y)|^2$. By Itô's formula, we have

$$d\xi_t = 2\sum_{i=1}^N \langle x_t(x) - x_t(y), A_i(x_t(x)) - A_i(x_t(y)) \rangle \, dw_t^i$$
$$+ 2\langle x_t(x) - x_t(y), \tilde{A}_0(x_t(x)) - \tilde{A}_0(x_t(y)) \rangle \, dt$$
$$+ \sum_{i=1}^N |A_i(x_t(x)) - A_i(x_t(y))|^2 \, dt.$$

The Itô stochastic contraction $d\xi_t \cdot d\xi_t$ is dominated by

$$4\sum_{i=1}^N \langle x_t(x) - x_t(y), A_i(x_t(x)) - A_i(x_t(y)) \rangle^2 \leq 4L_1^2 \xi_t^2.$$

Again, by Itô's formula,

$$d\xi_t^p = 2p\sum_{i=1}^N \xi_t^{p-1} \langle x_t(x) - x_t(y), A_i(x_t(x)) - A_i(x_t(y)) \rangle \, dw_t^i$$
$$+ 2p\xi_t^{p-1} \langle x_t(x) - x_t(y), \tilde{A}_0(x_t(x)) - \tilde{A}_0(x_t(y)) \rangle \, dt$$
$$+ p\xi_t^{p-1} \sum_{i=1}^N |A_i(x_t(x)) - A_i(x_t(y))|^2 \, dt + \frac{p(p-1)}{2} \xi_t^{p-2} \, d\xi_t \cdot d\xi_t,$$

which is less than

$$dM_t + (2pL_2 + 2p^2 L_1^2)\xi_t^p \, dt,$$

where $M_t$ is the martingale part of $\xi_t^p$. Taking expectations, we get

$$\mathbf{E}(\xi_t^p) \leq |x-y|^{2p} + (2pL_2 + 2p^2 L_1^2) \int_0^t \mathbf{E}(\xi_s^p) \, ds.$$

Now, Gronwall's lemma gives

$$\mathbf{E}(\xi_t^p) \leq |x-y|^{2p} e^{2pL_2 + 2p^2 L_1^2}, \qquad t \in [0,1],$$

which is nothing but (1.8). □

1.3. *Precise $L^p$-estimates for the two-point motion under local Lipschitz conditions.* We shall next assume that the vector fields $\tilde{A}_0, A_1, \ldots, A_N$ are only locally Lipschitz. We shall describe growth conditions in $m$ for the Lipschitz coefficients $L_m$ valid on Euclidean balls of radius $m$ that lead to $L^p$-moment estimates for the two-point motion of the flow. For this purpose, set

(1.9) $$L_{m,1}^2 = \sum_{i=1}^N \sup_{|x| \leq m} \|A_i'(x)\|^2, \qquad L_{m,2} = \sup_{|x| \leq m} \|\tilde{A}_0'(x)\|,$$

where $A'_i$ denotes the Jacobian of the mapping $x \to A_i(x)$. Then for any $x, y \in B(m) := \{z \in \mathbf{R}^d; |z| \leq m\}$, we have

$$\sum_{i=1}^{N} |A_i(x) - A_i(y)|^2 \leq L_{m,1}^2 |x-y|^2, \qquad |\tilde{A}_0(x) - \tilde{A}_0(y)| \leq L_{m,2} |x-y|.$$

Now, consider a family of smooth functions $\varphi_m : \mathbf{R}^d \to \mathbf{R}$ satisfying $0 \leq \varphi_m \leq 1$ and

(1.10)
$$\begin{aligned} \varphi_m(x) &= 1 \quad \text{for } |x| \leq m, \\ \varphi_m(x) &= 0 \quad \text{for } |x| > m+2, \\ \sup_m \sup_{x \in \mathbf{R}^d} |\varphi'_m(x)| &\leq 1. \end{aligned}$$

Define $A_{m,i} = \varphi_m A_i$ for $i = 1, \ldots, N$ and $A_{m,0} = \varphi_m \tilde{A}_0$. Then we have

(1.11) $$\sup_{x \in \mathbf{R}^d} |A'_{m,i}(x)|^2 \leq 2\left(\sup_{|x| \leq m+2} |A_i(x)|^2 + \sup_{|x| \leq m+2} \|A'_i(x)\|^2\right),$$

(1.12) $$\sup_{x \in \mathbf{R}^d} |A'_{m,0}(x)| \leq \sup_{|x| \leq m+2} |\tilde{A}_0(x)| + \sup_{|x| \leq m+2} \|\tilde{A}'_0(x)\|.$$

Set
$$\tilde{L}_{m,1}^2 = \sum_{i=1}^{N} \sup_{x \in \mathbf{R}^d} \|A'_{m,i}(x)\|^2, \qquad \tilde{L}_{m,2} = \sup_{x \in \mathbf{R}^d} \|A'_{m,0}(x)\|.$$

Let $(x_t^m(x))$ be the solution of the following stochastic differential equation:
$$dx_t^m = \sum_{i=1}^{N} A_{m,i}(x_t^m)\, dw_t^i + A_{m,0}(x_t^m)\, dt, \qquad x_0^m = x.$$

Applying (1.8), we get for $p \geq 2$,

(1.13) $$\mathbf{E}(|x_t^m(x) - x_t^m(y)|^{2p}) \leq |x-y|^{2p} e^{2p^2 \tilde{L}_{m,1}^2 + 2p\tilde{L}_{m,2}}, \qquad t \in [0,1].$$

We have
$$|x_t(x) - x_t(y)|^p = \sum_{m=1}^{+\infty} |x_t(x) - x_t(y)|^p \mathbf{1}_{\{m-1 \leq Y_1(x) \vee Y_1(y) < m\}}$$
$$= \sum_{m=1}^{+\infty} |x_t^m(x) - x_t^m(y)|^p \mathbf{1}_{\{m-1 \leq Y_1(x) \vee Y_1(y) < m\}}.$$

According to (1.13) and the Cauchy–Schwarz inequality, we obtain

(1.14)
$$\mathbf{E}(|x_t(x) - x_t(y)|^p)$$
$$\leq |x-y|^p \sum_{m=1}^{+\infty} e^{p^2 \tilde{L}_{m,1}^2 + p\tilde{L}_{m,2}} \sqrt{P(Y_1(x) \vee Y_1(y) \geq m-1)}.$$

With the aid of this inequality, we are able to formulate growth conditions on the Lipschitz constants ensuring global moment estimates for the flow. In the following theorems, this will be done consecutively under (H1) and (H2).

THEOREM 1.7. *Assume* (H1). *Let $p \geq 2$. Suppose that $L_{m,1} \leq \alpha m$, $L_{m,2} \leq \beta m^2$. For $R > 0$, let $\delta_0$ be given according to Corollary* 1.2. *Suppose that*

$$(1.15) \qquad p^2\alpha^2 + p\beta < \delta_0/2.$$

*Then for any $R > 0$, there exists a constant $C_{p,R} > 0$ such that*

$$(1.16) \quad \mathbf{E}(|x_t(x) - x_t(y)|^p) \leq C_{p,R}|x-y|^p \qquad \text{for } x, y \in B(R), t \in [0,1].$$

*In particular, if for some $\varepsilon > 0$ and constants $\beta_1, \beta_2$, we have*

$$L_{m,1} \leq \beta_1 m^{1-\varepsilon}, \qquad L_{m,2} \leq \beta_2 m^{2-\varepsilon},$$

*then for any $p \geq 2$, there exists $C_p > 0$ such that (1.16) holds.*

PROOF. Let $C_R = \sup_{|x| \leq R} \mathbf{E}(e^{\delta_0 Y_1^2(x)})$. Then for $m \geq 1$ and $x, y \in B(R)$,

$$\sqrt{P(Y_1(x) \vee Y_1(y) \geq m-1)} \leq \sqrt{2C_R} e^{-\delta_0(m-1)^2/2}.$$

On the other hand, by (1.11) and (1.12), we have

$$\tilde{L}_{m,1}^2 \leq 2NC_1^2 + 2\alpha^2(m+2)^2, \qquad \tilde{L}_{m,2} \leq \beta m^2 + (C_2 + 2\beta)m + 3C_2 + 4\beta.$$

Therefore, there exists a constant $\gamma_p > 0$, independent of $m$, such that

$$e^{p^2\tilde{L}_{m,1}^2 + p\tilde{L}_{m,2}} \leq \gamma_p e^{(p^2\alpha^2 + p\beta)m^2} e^{(2\alpha^2 + C_2 + 4\beta)m}.$$

Now, using (1.14), we get

$$\mathbf{E}(|x_t(x) - x_t(y)|^p)$$
$$\leq \gamma_p \sqrt{2C_R} |x-y|^p \sum_{m=1}^{+\infty} e^{-\delta_0(m-1)^2/2} e^{(p^2\alpha^2 + p\beta)m^2} e^{(2\alpha^2 + C_2 + 4\beta)m}.$$

It is clear that if $p^2\alpha^2 + p\beta < \delta_0/2$, then the above series converges, so that (1.16) follows. $\square$

REMARK. One can specify the $R$-dependence of the constant $C_{p,R}$ by looking at the proof of Corollary 1.2. It is seen that there is a subtle trade-off between $R$ and the parameter $\beta$ appearing in the bound for the Lipschitz constants $L_{m,2}$ which, in our setting, is expressed through the value of $\delta_0 = \delta_0(C_1, C_2, R)$.

Under (H2), the growth of the diffusion vector fields has to be counterbalanced by a slower growth of the local Lipschitz constants. We shall formulate them implicitly through conditions on the $\tilde{L}_{m,1}$ and $\tilde{L}_{m,2}$.

THEOREM 1.8. *Assume* (H2) *and the existence of constants* $\beta_1, \beta_2$ *such that*

(1.17) $$\tilde{L}^2_{m,1} \leq \beta_1 \log m, \qquad \tilde{L}_{m,2} \leq \beta_2 \log m.$$

*Then for any* $p \geq 2, R > 0$, *there exists a constant* $C_{p,R} > 0$ *such that*

(1.18) $$\mathbf{E}(|x_t(x) - x_t(y)|^p) \leq C_{p,R}|x - y|^p \qquad \text{for } x, y \in B(R), \ t \in [0,1].$$

PROOF. Let $q \geq 2$. By (1.3), $\alpha_{q,R} = \sup_{|x| \leq R} \mathbf{E}(Y_1(x)^q)$ is finite. Then for any $|x| \leq R$ and $m \geq 2$,

$$P(Y_1(x) \geq m - 1) \leq \alpha_{q,R} \frac{1}{(m-1)^q}.$$

On the other hand, under condition (1.17),

$$e^{p^2 \tilde{L}^2_{m,1} + p\tilde{L}_{m,2}} \leq (m+2)^{\beta_1 p^2 + p\beta_2}.$$

Therefore, if we take $\frac{q}{2} > \beta_1 p^2 + \beta_2 p + 2$, the series

$$\sum_{m \geq 2} \frac{1}{(m-1)^{q/2}} \cdot (m+2)^{\beta_1 p^2 + p\beta_2}$$

converges. Now, using (1.14), we obtain the desired result (1.18). □

**2. Moment estimates for regularized ordinary differential equations.** Let $n \geq 1$ be an integer. Define $(w^n_t)_{t \in [0,1]}$ by $w^n_0 = 0$ and

(2.1) $$\dot{w}^n_t = 2^n(w_{(\ell+1)2^{-n}} - w_{\ell 2^{-n}}) \qquad \text{for } t \in [\ell 2^{-n}, (\ell+1)2^{-n}[.$$

Let $x^n_t(x)$ be the solution of the following ordinary differential equation:

(2.2) $$dx^n_t = \sum_{i=1}^N A_i(x^n_t) \dot{w}^{n,i}_t \, dt + A_0(x^n_t) \, dt, \qquad x^n_0 = x.$$

The aim of this section is to prove moment estimates for one- and two-point motions of these regularized ordinary differential equations, uniformly in the discretization parameter $n$. For this purpose, we shall use the techniques presented in the previous section, involving the specification of Lipschitz constants.

2.1. *Uniform moment estimates for the one-point motions.* Define $Y_n(t,x) = \sup_{0 \leq s \leq t} |x^n_s(x)|$. Set

(2.3) $$B_{i,k} = \sum_{j=1}^d \frac{\partial A_i}{\partial x_j} A^j_k \qquad \text{for } i = 1, \ldots, N \text{ and } k = 0, 1, \ldots, N.$$

For the first uniform boundedness result, we shall work under growth assumptions very close to (H1) of the previous section.

PROPOSITION 2.1. *Assume that*

$$\sum_{i=1}^{N} |A_i(x)|^2 \leq C_1^2, \qquad |A_0(x)| \leq C_2(1+|x|) \tag{2.4}$$

*and*

$$|B_{ik}(x)| \leq C_3(1+|x|) \qquad \text{for all } i,k. \tag{2.5}$$

*Then there exist positive constants $\alpha_1$ and $\alpha_2$, independent of $n$ and $p$, such that*

$$\mathbf{E}(Y_n(1,x)^p) \leq (1+|x|)^p \alpha_1^p e^{\alpha_2 p^2}. \tag{2.6}$$

PROOF. For $t \in [0,1]$, define $t_n = k2^{-n}$ if $t \in [k2^{-n}, (k+1)2^{-n}[$ and $t_n^+ = t_n + 2^{-n}$. We then have, for fixed but arbitrary $t \in [0,1]$,

$$x_t^n = x + \sum_{i=1}^{N} \int_0^t A_i(x_{s_n}^n) \dot{w}_s^{n,i}\, ds + \int_0^t A_0(x_s^n)\, ds$$

$$+ \sum_{i=1}^{N} \int_0^t (A_i(x_s^n) - A_i(x_{s_n}^n)) \dot{w}_s^{n,i}\, ds$$

$$= x + M_n(t) + \int_0^t A_0(x_s^n)\, ds + R_n(t),$$

accordingly. Consider $Y_i(s) = A_i(x_{s_n}^n)$ for $s < t_n$ and $Y_i(s) = (t-t_n)2^n A_i(x_{t_n}^n)$ for $t_n \leq s \leq t$. Then $M_n(t) = \sum_{i=1}^{N} \int_0^{t_n^+} Y_i(s)\, dw_s^i$. We have

$$\int_0^{t_n^+} |Y_i(s)|^2\, ds = \int_0^{t_n} |Y_i(s)|^2\, ds + 2^{-n}(t-t_n)^2 2^{2n} |A_i(x_{t_n}^n)|^2 \leq \int_0^t |A_i(x_{s_n}^n)|^2\, ds$$

and by Burkholder's inequality,

(i) $\mathbf{E}(|M_n(t)|^p) \leq C\sqrt{p^p}\, \mathbf{E}\left[\left(\int_0^{t_n^+} \sum_{i=1}^{N} |Y_i(s)|^2\, ds\right)^{p/2}\right] \leq CC_1^p \sqrt{p^p}.$

Observe that for fixed $n$, $t \to M_n(t)$ is not a martingale. Only $k \to M_n(k2^{-n})$ is a $\mathcal{F}_{k2^{-n}}$-martingale. Let $t \in [\ell 2^{-n}, (\ell+1)2^{-n}[$. According to (i) and by Doob's maximal inequality, we have

(ii) $\mathbf{E}\left(\sup_{0 \leq k \leq \ell} |M_n(k2^{-n})|^p\right) \leq 2e\mathbf{E}(|M_n(t_n)|^p) \leq 2eCC_1^p\sqrt{p^p}.$

Here, $e$ is Euler's constant, resulting from the simple estimate

$$\left(\frac{p}{p-1}\right)^p \leq 2e, \qquad p > 1.$$

Now, for $s \in [k2^{-n}, (k+1)2^{-n}[$,

(iii) $M_n(s) = M_n(k2^{-n}) + (s - k2^{-n}) \sum_{i=1}^{N} A_i(x_{k2^{-n}}^n)(w_{(k+1)2^{-n}}^i - w_{k2^{-n}}^i) 2^n.$

Then $|M_n(s)| \leq |M_n(k2^{-n})| + C_1 2^{-n/2} \Gamma_n(k2^{-n})$, where

(2.7) $$\Gamma_n(s) = 2^{n/2} \sum_{i=1}^{N} |w_{s_n^+}^i - w_{s_n}^i|.$$

Therefore,

(iv) $\sup_{0 \leq s \leq t} |M_n(s)| \leq \sup_{0 \leq k \leq \ell} |M_n(k2^{-n})| + C_1 \sup_{0 \leq k \leq \ell} (2^{-n/2} \Gamma_n(k2^{-n})).$

Now using Lemma 2.2 below, we have for $p \geq 2$,

(v) $\mathbf{E}\left[\sup_{0 \leq k \leq \ell} (2^{-n/2} \Gamma(k2^{-n}))^p\right] \leq \sum_k 2^{-np/2} \mathbf{E}(\Gamma_n(k2^{-n})^p)$
$\leq 2^n \cdot 2^{-np/2} (CN)^p \sqrt{p^p} \leq (CN)^p \sqrt{p^p}.$

So, combining (iv), (ii) and (v), we finally obtain

(2.8) $$\left\| \sup_{0 \leq s \leq t} |M_n(s)|^p \right\|_p \leq CC_1 \sqrt{p}.$$

The remainder term $R_n$ is more delicate to estimate. Using the vector fields defined in (2.3), we may express $R_n$ by

$$R_n(t) = \sum_{i,k=1}^{N} \int_0^t \left[ \int_{s_n}^s B_{ik}(x_\sigma^n) \dot{w}_\sigma^{n,k} \dot{w}_s^{n,i} d\sigma \right] ds + \sum_{i=1}^{N} \int_0^t \left[ \int_{s_n}^s B_{i0}(x_\sigma^n) \dot{w}_s^{n,i} d\sigma \right] ds.$$

Let $R_{n,1}$ and $R_{n,2}$ be the two consecutive terms on the right-hand side of the preceding equation. Using hypothesis (2.4), for $\sigma \in [s_n, s[$, we obtain

$$|x_\sigma^n| \leq |x_{s_n}^n| + C_1 2^{-n} \sum_{i=1}^{N} |\dot{w}_{s_n}^{n,i}| + C_2 \int_{s_n}^\sigma (1 + |x_s^n|) ds.$$

Hence, Gronwall's lemma implies, with universal constants $C_1, C_2$, that

(2.9) $$1 + |x_\sigma^n| \leq (|x_{s_n}^n| + 1 + C_1 2^{-n/2} \Gamma_n(s_n)) e^{C_2 2^{-n}}.$$

Using (2.9) and hypothesis (2.5), we have

(2.10) $|R_{n,2}(t)| \leq C_3 e^{C_2 2^{-n}} \left[ \int_0^t (|x_{s_n}^n| + 1) \Gamma_n(s_n) ds + C_1 \int_0^t \Gamma_n(s_n)^2 ds \right].$

By independence of $x_{s_n}^n$ and $\Gamma_n(s_n)$, we have

(2.11) $\mathbf{E}((|x_{s_n}^n| + 1)^p \Gamma_n(s_n)^p) \leq \mathbf{E}((1 + Y_n(s,x))^p) \mathbf{E}(\Gamma_n(s_n)^p).$

Combining (2.10) and (2.11) and again using (2.14) in Lemma 2.2, we get

(2.12) $$\left\| \sup_{0 \leq s \leq t} |R_{n,2}(s)| \right\|_p$$
$$\leq C_3 e^{C_2} \left( CN\sqrt{p} \int_0^t (1 + \|Y_n(s,x)\|_p) ds + C_1 C^2 N^2 p \right).$$

In the same way,

$$
\left\| \sup_{0 \le s \le t} |R_{n,1}(s)| \right\|_p
$$
$$
\le C_3 e^{C_2} \left( C^2 N^2 p \int_0^t (1 + \|Y_n(s,x)\|_p) \, ds + C_1 C^3 N^3 p^{3/2} \right), \tag{2.13}
$$

where $C_3$ is another universal constant and $C$ results from Lemma 2.2. Now define $\psi(t) = \|Y_n(t,x)\|_p$. Combining (2.8), (2.12) and (2.13), we finally obtain

$$
\psi(t) + 1 \le |x| + 1 + CC_1\sqrt{p} + C_3 e^{C_2}(C_1 C^2 N^2 p + C_1 C^3 N^3 p^{3/2})
$$
$$
+ C_3 e^{C_2}(CN\sqrt{p} + C^2 N^2 p) \int_0^t (1 + \psi(s)) \, ds.
$$

From the structure of the bound just obtained, we see that there are two constants $\alpha_1, \alpha_2 > 0$, independent of $n$ and $p$, such that $\psi(1) \le (|x| + 1)\alpha_1 e^{\alpha_2 p}$ holds. The result (2.6) follows. □

LEMMA 2.2. *There exists a constant $C > 0$ such that*

$$
\|\Gamma_n(s)\|_q \le CN\sqrt{q} \quad \text{for all } s \in [0,1[,\ n \ge 1,\ q \ge 2. \tag{2.14}
$$

PROOF. Let $s \in [k2^{-n}, (k+1)2^{-n}[$ be given. Let $\gamma_i = 2^{n/2}(w^i_{(k+1)2^{-n}} - w^i_{k2^{-n}})$. Then $\gamma_1, \ldots, \gamma_N$ are independent standard Gaussian random variables. For any $1 \le i \le N$,

$$
\mathbf{E}(|\gamma_i|^q) = 2 \int_0^{+\infty} s^q e^{-s^2/2} \frac{ds}{\sqrt{2\pi}} = \frac{2^{q/2}}{\sqrt{\pi}} \int_0^{+\infty} s^{(q+1)/2 - 1} e^{-s} \, ds.
$$

By well-known properties of the Gamma function, the above quantity is dominated by $Cq^{q/2}$, with a universal constant $C > 0$. Now,

$$
\|\Gamma_n(s)\|_q \le \sum_{i=1}^N \|\gamma_i\|_{L^q} \le CN\sqrt{q}.
$$

We thus obtain (2.14). □

We next discuss the case where condition (2.4) is replaced by

$$
\sum_{i=1}^N |A_i(x)|^2 \le C_1^2(1 + |x|^2), \qquad |A_0(x)| \le C_2(1 + |x|). \tag{2.15}
$$

(2.15) combined with (2.5) resembles (H2) of the previous section.

PROPOSITION 2.3. *Assume* (2.15) *and* (2.5). *Then for any* $p \geq 2$, *there exists a constant* $C_p > 0$ *such that*

(2.16) $$\sup_{0 \leq t \leq 1} \mathbf{E}(|x_t^n(x)|^p) \leq C_p(1 + |x|^p) \quad \text{for any } n \geq 1.$$

PROOF. We resume the computation carried out in the proof of the previous proposition, taking into account the linear growth of coefficients $A_1, \ldots, A_N$. Let $t \in [0, 1]$ be fixed and set $M_n(t) = \sum_{i=1}^{N} \int_0^t A_i(x_{s_n}^n) \dot{w}_s^{n,i} \, ds$. By the previous computations, we see that for some constant $C_p > 0$,

(i) $\mathbf{E}(|M_n(t)|^p) \leq C_p \int_0^t (1 + \mathbf{E}(|x_{s_n}^n|^p)) \, ds.$

Moreover, for $\sigma \in [s_n, s_n^+[$, we have

$$|x_\sigma^n| \leq |x_{s_n}^n| + C_1\left(\int_{s_n}^\sigma (1 + |x_s^n|) \, ds\right) \sum_{i=1}^{N} |\dot{w}_s^{n,i}| + C_2 \int_{s_n}^\sigma (1 + |x_s^n|) \, ds.$$

So Gronwall's lemma gives, with some universal constants $C_1, C_2$,

$$|x_\sigma^n| + 1 \leq (|x_{s_n}^n| + 1) e^{2^{-n}(C_2 + C_1 \sum_{i=1}^{N} |\dot{w}_{s_n}^{n,i}|)}.$$

It follows that

(2.17) $$|x_\sigma^n| + 1 \leq e^{C_2}(|x_{s_n}^n| + 1) e^{C_1 \Gamma_n(s_n)}, \quad \sigma \in [s_n, s_n^+[.$$

Replacing (2.9) by (2.17) in the estimate of $R_n(t)$, we have, with another universal constant $C_3$,

$$|R_{n,2}(t)| \leq C_3 e^{C_2} \int_0^t (|x_{s_n}^n| + 1) e^{C_1 \Gamma_n(s_n)} \Gamma_n(s_n) \, ds,$$

$$|R_{n,1}(t)| \leq C_3 e^{C_2} \int_0^t (|x_{s_n}^n| + 1) e^{C_1 \Gamma_n(s_n)} \Gamma_n(s_n)^2 \, ds.$$

By a direct calculation,

(ii) $\mathbf{E}(e^{2pC_1 \Gamma_n(s)}) \leq 2^N e^{4p^2 C_1^2 N/2}.$

Now, using the independence of $|x_{s_n}^n|$ and $\Gamma(s_n)$, (ii) and (2.14), we see that there exists a constant $C_p > 0$ such that

(iii) $\mathbf{E}(|R_n(t)|^p) \leq C_p \int_0^t (1 + \mathbf{E}(|x_{s_n}^n|^p)) \, ds.$

Therefore, (i) and (iii) imply that

$$\mathbf{E}(|x_t^n|^p) \leq C_p\left(|x|^p + \int_0^t (1 + \mathbf{E}(|x_{s_n}^n|^p)) \, ds + \int_0^t (1 + \mathbf{E}(|x_s^n|^p)) \, ds\right).$$

Finally, we consider $\psi(t) = \sup_{0 \leq s \leq t} \mathbf{E}(|x_s^n|^p) + 1$. The inequality just derived implies that

$$\psi(t) \leq C_p(|x|^p + 1) + 2C_p \int_0^t \psi(s) \, ds.$$

So, a final application of Gronwall's lemma yields another constant $C_p$ such that

$$\sup_{0 \leq t \leq 1} \mathbf{E}(|x_t^n|^p) \leq C_p(1 + |x|^p). \qquad \square$$

Using the same techniques, we may also derive uniform moment estimates for the time fluctuations of the approximate ordinary differential equations.

PROPOSITION 2.4. *Assume* (2.15) *and* (2.5). *Then for any* $p \geq 2$, *there exists a constant* $C_p > 0$, *independent of* $n$, *such that*

$$(2.18) \qquad \mathbf{E}(|x_s^n(x) - x_t^n(x)|^p) \leq C_p(1 + |x|^p)|s - t|^{p/2}.$$

Finally, we derive a result describing a bound for the maximal growth of the one-point motions of the regularizing ordinary differential equations, uniformly in $n$.

THEOREM 2.5. *Assume* (2.15) *and* (2.5). *Then for any* $p \geq 2$, *there exists a constant* $C_p > 0$ *such that*

$$(2.19) \qquad \mathbf{E}(Y_n(1,x)^p) \leq C_p(1 + |x|^p) \qquad \text{for any } n \geq 1.$$

PROOF. Let $\gamma > 0$ be a parameter such that $0 < \gamma < 1/2$ and $q \geq 2$ be an integer such that $2q\gamma > 1, 2q(\frac{1}{2} - \gamma) > 1$. Then it is known from the regularity lemma of Garsia, Rodemich and Rumsey that

$$\sup_{0 \leq t \leq 1} |\psi(t)|^{2q} \leq C_{q,\gamma} \int_0^1 \int_0^1 \frac{|\psi(s) - \psi(t)|^{2q}}{|t - s|^{1+2q\gamma}} \, ds \, dt.$$

Therefore, we have

$$\mathbf{E}\left(\sup_{0 \leq t \leq 1} |x_t^n(x)|^{2qp}\right) \leq C_{q,\gamma}^p \int_0^1 \int_0^1 \frac{\mathbf{E}(|x_s^n(x) - x_t^n(x)|^{2qp})}{|t - s|^{(1+2q\gamma)p}} \, ds \, dt.$$

But, by (2.18), this bound is dominated by $C_p(1 + |x|^p)^{2q}$ since

$$\int_0^1 \int_0^1 |t - s|^{qp - (1+2q\gamma)p} \, ds \, dt \leq 1.$$

So we get (2.19). $\square$

2.2. *Uniform moment estimates for the two-point motions.* For vector fields satisfying global Lipschitz conditions and regularizations as considered here, Bismut [1] or Moulinier [10] proved that $\mathbf{E}(|x_t^n(x) - x_t^n(y)|^p) \leq C_p|x - y|^p$ for all $x, y \in \mathbf{R}^d$, where $C_p$ is independent of $n$. However, the dependence of $C_p$ on the Lipschitz continuity properties of the vector fields is not specified. In what follows, we shall make this functional dependence explicit.

THEOREM 2.6. *Assume that for $x, y \in \mathbf{R}^d$,*

$$(2.20) \quad \sum_{i=1}^{N} |A_i(x) - A_i(y)|^2 \leq L_1^2 |x-y|^2, \qquad |A_0(x) - A_0(y)| \leq L_2 |x-y|$$

*and for all $1 \leq i \leq N, 1 \leq k \leq N$,*

$$(2.21) \quad |B_{ik}(x) - B_{ik}(y)| \leq K_1 |x-y|, \qquad |B_{i0}(x) - B_{i0}(y)| \leq K_2 |x-y|.$$

*Let $C$ be the constant appearing in Lemma 2.2. Define*

$$\alpha_n = 2p((2p-1)L_1^2 + K_1)(4C^2 N^2 2^N e^{8p^2 N 2^{-n} L_1^2}) e^{2p 2^{-n} L_2}$$
$$+ 2^{-n/2} 2p((2p-1)L_1 L_2 + K_2)(2CN 2^N e^{8p^2 N 2^{-n} L_1^2}) e^{2p 2^{-n} L_2}.$$

*Then*

$$\mathbf{E}(|x_t^n(x) - x_t^n(y)|^{2p}) \leq |x-y|^{2p} e^{2pL_2} e^{\alpha_n} \leq |x-y|^{2p} e^{2pL_2} e^{\alpha_1}.$$

PROOF. For $n, x, y, t$ fixed, we have

$$x_t^n(x) - x_t^n(y) = x - y + \sum_{i=1}^{N} \int_0^t (A_i(x_s^n(x)) - A_i(x_s^n(y))) \dot{w}_s^{n,i} \, ds$$
$$+ \int_0^t (A_0(x_s^n(x)) - A_0(x_s^n(y))) \, ds.$$

Set $\xi_t = |x_t^n(x) - \xi_t^n(y)|^2$. Then

$$d\xi_t = 2 \sum_{i=1}^{N} \langle x_t^n(x) - x_t^n(y), A_i(x_t^n(x)) - A_i(x_t^n(y)) \rangle \dot{w}_t^{n,i} \, dt$$
$$+ 2 \langle x_t^n(x) - x_t^n(y), \ A_0(x_t^n(x)) - A_0(x_t^n(y)) \rangle \, dt.$$

Set

$$Q_i(t) = \langle x_t^n(x) - x_t^n(y), A_i(x_t^n(x)) - A_i(x_t^n(y)) \rangle \qquad \text{for } i = 0, 1, \ldots, N.$$

Then $d\xi_t$ has the decomposition $d\xi_t = 2 \sum_{i=1}^{N} Q_i(t) \dot{w}_t^{n,i} \, dt + 2Q_0(t) \, dt$. For $p \geq 2$, we have

$$d\xi_t^p = 2p \sum_{i=1}^{N} \xi_t^{p-1} Q_i(t) \dot{w}_t^{n,i} \, dt + 2p \xi_t^{p-1} Q_0(t) \, dt$$

$$(2.22) \qquad = 2p \sum_{i=1}^{N} \xi_{t_n}^{p-1} Q_i(t_n) \dot{w}_t^{n,i} \, dt + 2p \xi_t^{p-1} Q_0(t) \, dt$$

$$+ 2p \sum_{i=1}^{N} (\xi_t^{p-1} Q_i(t) - \xi_{t_n}^{p-1} Q_i(t_n)) \dot{w}_t^{n,i} \, dt.$$

Let $M_t = 2p \sum_{i=1}^N \int_0^t \xi_{s_n}^{p-1} Q_i(s_n) \dot{w}_s^{n,i} \, ds$. Then $\mathbf{E}(M_t) = 0$. Moreover, we have

$$(2.23) \qquad 2p \int_0^t |\xi_s^{p-1} Q_0(s)| \, ds \leq 2pL_2 \int_0^t \xi_s^p \, ds.$$

To estimate the third term $R(t) = 2p \sum_{i=1}^N \int_0^t (\xi_s^{p-1} Q_i(s) - \xi_{s_n}^{p-1} Q_i(s_n)) \dot{w}_s^{n,i} \, ds$ appearing on the right-hand side of (2.22), we compute the derivative of $\xi_s^{p-1} Q_i(s)$. We get

$$(\xi_s^{p-1} Q_i(s))' = (p-1) \xi_s^{p-2} \xi_s' Q_i(s) + \xi_s^{p-1} Q_i'(s).$$

Computing $Q_i'(s)$ and using our Lipschitz continuity hypotheses, we get

$$|Q_i'(s)| \leq (K_1 + L_1^2) \xi_s \sum_{k=1}^N |\dot{w}_s^{n,k}| + (L_1 L_2 + K_2) \xi_s.$$

Therefore,

$$(2.24) \quad |(\xi_s^{p-1} Q_i(s))'|$$
$$\leq ((2p-1)L_1 L_2 + K_2) \xi_s^p + ((2p-1)L_1^2 + K_1) \xi_s^p \sum_{k=1}^N |\dot{w}_s^{n,k}|.$$

To estimate the contribution of $\xi_s^p$, first note that for $\sigma \in [s_n, s_n^+[$, we have

$$|x_\sigma^n(x) - x_\sigma^n(y)| \leq |x_{s_n}^n(x) - x_{s_n}^n(y)| + L_1 \left( \int_{s_n}^\sigma |x_u^n(x) - x_u^n(y)| \, du \right) \sum_{i=1}^N |\dot{w}_{s_n}^{n,i}|$$
$$+ L_2 \int_{s_n}^\sigma |x_u^n(x) - x_u^n(y)| \, du.$$

Now apply Gronwall's lemma. This leads to

$$|x_\sigma^n(x) - x_\sigma^n(y)| \leq |x_{s_n}^n(x) - x_{s_n}^n(y)| \cdot e^{2^{-n} L_1 \sum_{i=1}^N |\dot{w}_{s_n}^{n,i}|} e^{2^{-n} L_2}.$$

Therefore, for $\sigma \in [s_n, s_n^+[$,

$$(2.25) \qquad \xi_\sigma^p \leq \xi_{s_n}^p \cdot e^{2p 2^{-n/2} L_1 \Gamma_n(s_n)} e^{2p 2^{-n} L_2}.$$

Hence, by (2.24),

$$|R(t)| \leq 2p \sum_{i=1}^N \int_0^t \int_{s_n}^s |(\xi_\sigma^{p-1} Q_i(\sigma))'| |\dot{w}_s^{n,i}| \, d\sigma \, ds$$

$$\leq 2p \left\{ ((2p-1)L_1^2 + K_1) \int_0^t \int_{s_n}^s \xi_\sigma^p \left( \sum_{i=1}^N |\dot{w}_\sigma^{n,i}| \right) \left( \sum_{k=1}^N |\dot{w}_s^{n,k}| \right) d\sigma \, ds \right.$$

$$\left. + ((2p-1)L_1 L_2 + K_2) \int_0^t \int_{s_n}^s \xi_\sigma^p \left( \sum_{i=1}^N |\dot{w}_\sigma^{n,i}| \right) d\sigma \, ds \right\}$$

which, according to (2.25), is dominated by

$$2p\Big\{((2p-1)L_1^2 + K_1)e^{2p2^{-n}L_2}\int_0^t \xi_{s_n}^p \Gamma_n(s_n)^2 e^{2p2^{-n/2}L_1\Gamma_n(s_n)}\,ds$$
$$+ 2^{-n/2}((2p-1)L_1L_2 + K_2)e^{2p2^{-n}L_2}\int_0^t \xi_{s_n}^p \Gamma_n(s_n) e^{2p2^{-n/2}L_1\Gamma_n(s_n)}\,ds\Big\}.$$

Next, we employ the independence of $\xi_{s_n}$ and $\Gamma_n(s_n)$. Therefore,

$$\mathbf{E}(\xi_{s_n}^p \Gamma_n(s_n)^2 e^{2p2^{-n/2}L_1\Gamma_n(s_n)}) = \mathbf{E}(\xi_{s_n}^p)\mathbf{E}(\Gamma_n(s_n)^2 e^{2p2^{-n/2}L_1\Gamma_n(s_n)}).$$

By estimates derived above, using Lemma 2.2, we have

$$\mathbf{E}(\Gamma_n(s_n)^2 e^{2p2^{-n/2}L_1\Gamma_n(s_n)}) \leq 4C^2 N^2 2^N e^{8p^2 N 2^{-n}L_1^2}$$

and

$$\mathbf{E}(\Gamma_n(s_n) e^{2p2^{-n/2}L_1\Gamma_n(s_n)}) \leq 2CN 2^N e^{8p^2 N 2^{-n}L_1^2}.$$

Summarizing, the definition

$$(2.26) \quad \begin{aligned}\alpha_n &= 2p((2p-1)L_1^2 + K_1)(4C^2 N^2 2^N e^{8p^2 N 2^{-n}L_1^2})e^{2p2^{-n}L_2} \\ &\quad + 2^{-n/2}2p((2p-1)L_1L_2 + K_2)(2CN 2^N e^{8p^2 N 2^{-n}L_1^2})e^{2p2^{-n}L_2}\end{aligned}$$

implies the following inequality for $\mathbf{E}(|R(t)|)$:

$$\mathbf{E}(|R(t)|) \leq \alpha_n \int_0^t \mathbf{E}(\xi_{s_n}^p)\,ds.$$

Substituting all the estimates obtained so far in (2.22), we obtain

$$\mathbf{E}(\xi_t^p) \leq |x-y|^{2p} + 2pL_2 \int_0^t \mathbf{E}(\xi_s^p)\,ds + \alpha_n \int_0^t \mathbf{E}(\xi_{s_n}^p)\,ds.$$

Finally, let $\psi_u = \sup_{0\leq s\leq u}\mathbf{E}(\xi_s^p)$. For $T > 0$ and any $0 \leq t \leq T$, the above inequality then leads to

$$\mathbf{E}(\xi_t^p) \leq |x-y|^{2p} + 2pL_2\int_0^T \psi_s\,ds + \alpha_n\int_0^t \psi_s\,ds,$$

also expressible as $\psi_T \leq |x-y|^{2p} + (2pL_2 + \alpha_n)\int_0^T \psi_s\,ds$. So Gronwall's lemma implies that for any $0 \leq t \leq 1$,

$$\mathbf{E}(\xi_t^p) \leq |x-y|^{2p} e^{2pL_2} e^{\alpha_n}.$$

We thus have the desired result. $\square$

**3. Limit theorem without global Lipschitz conditions.** The expression (2.26) for $\alpha_n$ is quite complicated. But it gives the explicit dependence of our uniform moment estimates on the Lipschitz constants for the vector fields of the underlying stochastic differential equation. We shall exploit this fact in the present section to derive a theorem about the convergence of the ordinary differential equation regularizations given in the preceding section to the solution of the stochastic differential equation. The explicit form of the dependence allows us to relax the global Lipschitz conditions to suitable local ones. For this purpose, the techniques explained in the first section will be applied. Let us first formulate convenient local Lipschitz conditions.

Let $A_1, \ldots, A_N$ be $\mathcal{C}^2$-vector fields on $\mathbf{R}^d$ and $A_0$ be a $\mathcal{C}^1$-vector field. Suppose for $x, y \in B(n)$,

$$
\begin{aligned}
\sum_{i=1}^{N} |A_i(x) - A_i(y)|^2 &\leq L_{n,1}^2 |x-y|^2, \\
|A_0(x) - A_0(y)| &\leq L_{n,2} |x-y|,
\end{aligned}
\tag{3.1}
$$

with positive constants $L_{n,1}, L_{n,2}$. Choose a family of smooth functions $\varphi_n : \mathbf{R}^d \to \mathbf{R}$ satisfying $0 \leq \varphi_n \leq 1$ and

$$
\begin{aligned}
\varphi_n &= 1 \quad \text{on } B(n), \qquad \varphi_n = 0 \quad \text{on } B(n+2)^c, \\
\sup_n \|\varphi_n'\|_\infty &\leq 1, \qquad \sup_n \|\varphi_n''\|_\infty \leq C < +\infty,
\end{aligned}
\tag{3.2}
$$

where $\|\cdot\|_\infty$ denotes the uniform norm. Introduce the vector fields

$$A_{n,i} = \varphi_n A_i \qquad \text{for } i = 0, 1, \ldots, N.$$

Put

$$
\tilde{L}_{n,1}^2 = \sum_{i=1}^{N} \sup_{x \in \mathbf{R}^d} \|A'_{n,i}(x)\|^2, \qquad \tilde{L}_{n,2} = \sup_{x \in \mathbf{R}^d} \|A'_{n,0}(x)\|.
\tag{3.3}
$$

Define

$$B_{ik}^n = \sum_{j=1}^{d} \frac{\partial A_{n,i}}{\partial x_j} A_{n,k}^j \qquad \text{for } i = 1, \ldots, N \text{ and } k = 0, 1, \ldots, N$$

and set

$$
K_{n,1} = \sup_{i,k} \sup_{x \in \mathbf{R}^d} \|(B_{ik}^n)'(x)\|, \qquad K_{n,2} = \sup_{i} \sup_{x \in \mathbf{R}^d} \|(B_{i0}^n)'(x)\|.
\tag{3.4}
$$

For $n \in \mathbf{N}$, let $(z_t^n(x))$ be the solution of the ordinary differential equation

$$
dz_t^n = \sum_{i=1}^{N} A_{n,i}(z_t^n) \dot{w}_t^{n,i} \, dt + A_{n,0}(z_t^n) \, dt, \qquad z_0^n = x,
\tag{3.5}
$$

with $\dot{w}^{n,i}$ as defined in (2.1). We can apply Theorem 2.6 to obtain the estimate

(3.6) $$\mathbf{E}(|z_t^n(x) - z_t^n(y)|^{2p}) \leq |x-y|^{2p} e^{2p\tilde{L}_{n,2}} e^{\tilde{\alpha}_n},$$

where

(3.7) $$\begin{aligned}\tilde{\alpha}_n &= 2p((2p-1)\tilde{L}_{n,1}^2 + K_{n,1})(4C^2 N^2 2^N e^{8p^2 N 2^{-n}\tilde{L}_{n,1}^2}) e^{2p2^{-n}\tilde{L}_{n,2}} \\ &\quad + 2^{-n/2} 2p((2p-1)\tilde{L}_{n,1}\tilde{L}_{n,2} + K_{n,2}) \\ &\quad \times (2CN2^N e^{8p^2 N 2^{-n}\tilde{L}_{n,1}^2}) e^{2p2^{-n}\tilde{L}_{n,2}}.\end{aligned}$$

Now, suppose that with positive constants $\tilde{\beta}_i, 1 \leq i \leq 4$, we have

(3.8) $$\begin{aligned}\tilde{L}_{n,1}^2 &\leq \tilde{\beta}_1 \log n, & \tilde{L}_{n,2} &\leq \tilde{\beta}_2 \log n, \\ K_{n,1} &\leq \tilde{\beta}_3 \log n, & K_{n,2} &\leq \tilde{\beta}_4 (\log n)^{3/2}.\end{aligned}$$

Under these conditions, it is easy to see from the definition of $\tilde{\alpha}_n$ that there exists a constant $C_p$, independent of $n$, such that

(3.9) $$\tilde{\alpha}_n \leq C_p(\tilde{L}_{n,1}^2 + K_{n,1} + 1).$$

Therefore, (3.6) implies that

$$\mathbf{E}(|z_t^n(x) - z_t^n(y)|^{2p}) \leq |x-y|^{2p} e^{C_p} e^{2p\tilde{L}_{n,2}} e^{C_p(\tilde{L}_{n,1}^2 + K_{n,1})}.$$

Our aim is to obtain an estimate which is uniform relative to $n$. For this purpose, we shall again use the cut-off functions $\varphi_m$ introduced in (3.2). For the sake of simplicity, we shall formulate conditions only on the coefficients $A_0, A_1, \ldots, A_N$. For $m \geq 1$, set

$$C_{m,1}^2 = \sum_{i=1}^N \left(\sup_{|x| \leq m} |A_i(x)|^2\right), \qquad C_{m,2} = \sup_{|x| \leq m} |A_0(x)|,$$

$$J_{m,1} = \sup_{i,k \neq 0} \left(\sup_{|x| \leq m} \|B'_{ik}(x)\|^2\right), \qquad J_{m,2} = \sup_i \sup_{|x| \leq m} \|B'_{i0}(x)\|.$$

We shall work under the following hypotheses:

(H) $$\begin{cases} C_{m,1}^2 \leq \gamma_1 \log m, & C_{m,2} \leq \gamma_2 \log m, \\ L_{m,1}^2 \leq \beta_1 \log m, & L_{m,2} \leq \beta_2 \log m, \\ J_{m,1} \leq \delta_1 \log m, & J_{m,2} \leq \delta_2 (\log m)^{3/2}. \end{cases}$$

Recall that $A_{n,i} = \varphi_n A_i$. Under hypothesis (H), we have

$$\sum_{i=1}^N |A_{n,i}|^2 \leq \gamma_1 \log(n+2), \qquad |A_{n,0}| \leq \gamma_2 \log(n+2),$$

$$\sum_{i=1}^N \|A'_{n,i}\|^2 \leq 2(\gamma_1 + \beta_1)\log(n+2), \qquad \|A'_{n,0}\| \leq (\gamma_2 + \beta_2)\log(n+2).$$

Since $B^n_{ik} = \sum_{j=1}^d \frac{\partial \varphi_n}{\partial x_j} \varphi_n A_i A^j_k + \varphi_n^2 B_{ik}$, hypothesis (H) moreover implies that

$$\|(B^n_{ik})'\| \leq \tilde{\delta}_1 \log(n+2), \qquad \|(B^n_{i0})'\| \leq \tilde{\delta}_2 (\log(n+2))^{3/2},$$

for some constants $\tilde{\delta}_1$ and $\tilde{\delta}_2$. Therefore, hypothesis (H) implies conditions (3.8), so that (3.9) is validated. Now let $m \geq 1$. Consider

$$A_{m,n,i} = \varphi_m A_{n,i} \qquad \text{for } i = 0, 1, \ldots, N.$$

We have

$$(3.10) \quad \sum_{i=1}^N |A_{m,n,i}|^2 \leq \gamma_1 \log(m \wedge n + 2), \qquad |A_{m,n,0}| \leq \gamma_2 \log(m \wedge n + 2),$$

$$(3.11) \quad \sum_{i=1}^N \|A'_{m,n,i}\|^2 \leq \tilde{\beta}_1 \log(m \wedge n + 2), \qquad \|A'_{m,n,0}\| \leq \tilde{\beta}_1 \log(m \wedge n + 2)$$

and

$$(3.12) \quad \|(B^{mn}_{ik})'\| \leq \tilde{\delta}_1 \log(m \wedge n + 2), \qquad \|(B^{mn}_{i0})'\| \leq \tilde{\delta}_2 (\log(m \wedge n + 2))^{3/2}.$$

Let $(z^{mn}_t(x))$ be the solution of

$$(3.13) \quad dz^{mn}_t = \sum_{i=1}^N A_{m,n,i}(z^{mn}_t) \dot{w}^{n,i}_t \, dt + A_{m,n,0}(z^{mn}_t) \, dt, \qquad z^{mn}_0 = x.$$

Using (3.10)–(3.12) to estimate $\tilde{\alpha}_m$ in (3.7), we have for $m \leq n$,

$$\tilde{\alpha}_m \leq C_p((\tilde{\beta}_1 + \tilde{\gamma}_1) \log(m+2) + 1).$$

We conclude that

$$(3.14) \quad \begin{aligned} \mathbf{E}(|z^{mn}_t(x) - z^{mn}_t(y)|^{2p}) &\leq e^{C_p} e^{2p\tilde{\beta}_2 \log(m+2)} e^{C_p(\tilde{\beta}_1 + \tilde{\delta}_1) \log(m+2)} |x-y|^{2p} \\ &= e^{C_p}(m+2)^{2p\tilde{\beta}_2 + C_p(\tilde{\beta}_1 + \tilde{\delta}_1)} |x-y|^{2p}. \end{aligned}$$

Extrapolating in $m$ by means of the techniques presented in Section 1, we obtain the following moment estimate for the two-point motion, uniformly in the regularization parameter:

THEOREM 3.1. *Under the hypothesis* (H), *for any $p \geq 2$ and $R > 0$, there exists a constant $C_{p,R} > 0$, independent of $n$, such that*

$$(3.15) \qquad \mathbf{E}(|z^n_t(x) - z^n_t(y)|^p) \leq C_{p,R} \, |x-y|^p \qquad \text{for } x, y \in B(R).$$

PROOF. It is clear that (H) implies the growth conditions (2.15) and (2.5). Let $Y_n(x) = \sup_{0 \leq t \leq 1} |z_t^n(x)|$. We have

$$|z_t^n(x) - z_t^n(y)|^p = \sum_{m \geq 1} |z_t^n(x) - z_t^n(y)|^p \mathbf{1}_{\{m-1 \leq Y_n(x) \vee Y_n(y) < m\}}$$

$$= \sum_{m \geq 1} |z_t^{mn}(x) - z_t^{mn}(y)|^p \mathbf{1}_{\{m-1 \leq Y_n(x) \vee Y_n(y) < m\}}.$$

Let $q \geq 2$. By (2.16), there is a constant $C_{q,R} > 0$ such that for all $|x| \leq R$ and $|y| \leq R$,

$$P(Y_n(x) \vee Y_n(y) \geq m - 1) \leq C_{q,R} \frac{1}{m^q}.$$

Using (3.14), we have

$$\mathbf{E}(|z_t^{mn}(x) - z_t^{mn}(y)|^p \mathbf{1}_{\{m-1 \leq Y_n(x) \vee Y_n(y) < m\}})$$

$$\leq e^{C_p}(m+2)^{p\tilde{\beta}_2 + C_p(\tilde{\beta}_1 + \tilde{\delta}_1)/2} \cdot \sqrt{C_{q,R}} \frac{1}{m^{q/2}} |x-y|^p.$$

Now taking $q/2 \geq p\tilde{\beta}_2 + \frac{1}{2}C_p(\tilde{\beta}_1 + \tilde{\delta}_1) + 2$ gives (3.15). □

The following proposition states a similar uniform moment estimate for the time fluctuations of the solutions of the regularized equations:

PROPOSITION 3.2. *Assume hypothesis* (H) *is satisfied. For any* $p \geq 2$ *and* $R > 0$, *there exists a constant* $C_{p,R} > 0$, *independent of* $n$, *such that*

(3.16) $\mathbf{E}(|z_t^n(x) - z_s^n(x)|^p) \leq C_{p,R} |t-s|^{p/2}, \qquad |x| \leq R, s,t \in [0,1].$

PROOF. The coefficients $A_{n,i}$ and $B_{ik}^n$ satisfy (2.5) and (2.15). Therefore, we can apply Proposition 2.4 to obtain (3.16). □

We are finally in a position to prove the convergence of the ordinary differential equations' regularizations $(z_t^n)$ to the solution of the stochastic differential equation $(x_t)$ in the $L^p$ sense, uniformly in space and time. To state this result, we first establish it in a weaker sense.

LEMMA 3.3. *Let* $R > 0$ *and* $p \geq 2$. *Then*

(3.17) $$\lim_{n \to +\infty} \sup_{|x| \leq R} \sup_{0 \leq t \leq 1} \mathbf{E}(|z_t^n(x) - x_t(x)|^p) = 0.$$

PROOF. Let $Y_n(x) = \sup_{0 \leq t \leq 1} |z_t^n(x)|$ and $Y(x) = \sup_{0 \leq t \leq 1} |x_t(x)|$. Let $m \geq 1$. We have

$$\mathbf{E}(|z_t^n(x) - x_t(x)|^p) = \mathbf{E}(|z_t^n(x) - x_t(x)|^p \mathbf{1}_{\{Y_n(x) \vee Y(x) \leq m\}})$$
$$+ \mathbf{E}(|z_t^n(x) - x_t(x)|^p \mathbf{1}_{\{Y_n(x) \vee Y(x) > m\}}).$$

Due to (1.4) and (2.16), the second term is majorized by

$$C_p\mathbf{E}((Y_n(x)^p + Y(x)^p)\mathbf{1}_{\{Y_n(x)\vee Y(x)>m\}}) \leq C_{p,R}\frac{1}{\sqrt{m}}.$$

To obtain (3.17), it is therefore sufficient to prove that

$$(3.18) \quad \lim_{n\to+\infty}\sup_{|x|\leq R}\sup_{0\leq t\leq 1}\mathbf{E}(|z_t^n(x) - x_t(x)|^p\mathbf{1}_{\{Y_n(x)\vee Y(x)\leq m\}}) = 0.$$

Let $n > m + 2$. By uniqueness of solutions, on the subset $\{w;\ Y_n(x) \leq m\}$, $z_t^n(x) = x_t^{mn}(x)$ for all $t \in [0,1]$, where $x_t^{nm}(x)$ is the solution of the following ordinary differential equation:

$$dx_t^{nm} = \sum_{i=1}^{N}(\varphi_m A_i)(x_t^{nm}(x))\dot{w}_t^{n,i}\,dt + (\varphi_m A_0)(x_t^{nm}(x))\,dt, \qquad x_0^{nm} = x.$$

On the other hand, let $\tau_m(x) = \inf\{t > 0,\ |x_t(x)| \geq m\}$. Then $x_{t\wedge\tau_m(x)}(x)$ satisfies the following Itô stochastic differential equation:

$$dx_t^m(x) = \sum_{i=1}^{N}(\varphi_m A_i)(x_t^m(x))\,dw_t^i$$
$$+ \left(\varphi_m A_0 + \frac{1}{2}\sum_{i,j}\frac{\partial(\varphi_m A_i)}{\partial x_j}(\varphi_m A_i^j)\right)dt, \qquad x_0^m(x) = x.$$

It follows that on the subset $\{Y(x) \leq m\}$ or $\{\tau_m(x) \geq 1\}$, we have $x_t^m(x) = x_t(x)$ for all $t \in [0,1]$. Therefore,

$$\mathbf{E}(|z_t^n(x) - x_t(x)|^p\mathbf{1}_{\{Y_n(x)\vee Y(x)\leq m\}}) = \mathbf{E}(|x_t^{nm}(x) - x_t^m(x)|^p\mathbf{1}_{\{Y_n(x)\vee Y(x)\leq m\}})$$
$$\leq \mathbf{E}(|x_t^{nm}(x) - x_t^m(x)|^p).$$

We are now in the classical situation. Therefore, Moulinier's result from [10] can be applied to obtain (3.18). The proof of (3.17) is thus complete. $\square$

We finally strengthen the previous result to moment convergence, uniformly in space and time.

THEOREM 3.4. *Assume hypothesis* (H). *Then any* $p \geq 2$,

$$(3.19) \quad \lim_{n\to+\infty}\mathbf{E}\left(\sup_{0\leq t\leq 1}\sup_{|x|\leq R}|z_t^n(x) - x_t(x)|^p\right) = 0.$$

PROOF. Let $p \geq 2$ be given. By (3.15), (3.16) and the Kolmogorov modification theorem, there exists $\beta > 0$ such that for $|x| \leq R, |y| \leq R$ and $t, s \in [0,1]$,

$$(3.20) \quad |z_t^n(x) - z_s^n(y)| \leq F_n \cdot (|x-y|^\beta + |t-s|^\beta), \qquad n \geq 1,$$

where $\{F_n; n \geq 1\}$ is a family of measurable functions bounded in $L^p$ for any $p$. In the same way, according to Corollary 1.4 and Proposition 1.6, there exists $F \in L^p$ such that

(3.21) $$|x_t(x) - x_s(y)| \leq F \cdot (|x-y|^\beta + |t-s|^\beta).$$

Let $\varepsilon_n = \sup_{0 \leq t \leq 1} \sup_{|x| \leq R} \mathbf{E}(|z_t^n(x) - x_t(x)|^p)$. By Lemma 3.3, $\lim_{n \to +\infty} \varepsilon_n = 0$. Let $\sigma_n > 0$. Then there exist $N_n \leq C \, (\frac{1}{\sigma_n})^{d+1}$ points $x_1, \ldots, x_{N_n}$ in the ball $B(R)$ and $t_1, \ldots, t_{N_n} \in [0,1]$ such that

$$[0,1] \times B(R) \subset \bigcup_{i=1}^{N_n} [t_i - \sigma_n, t_i + \sigma_n] \times \{x; \ |x - x_i| \leq \sigma_n\}.$$

Let $(t, x) \in [0, 1] \times B(R)$. There exists one $i$ such that $|t - t_i| \leq \sigma_n$ and $|x - x_i| \leq \sigma_n$. We have, according to (3.20) and (3.21),

$$|z_t^n(x) - x_t(x)| \leq |z_t^n(x) - z_{t_i}^n(x_i)| + |z_{t_i}^n(x_i) - x_{t_i}(x_i)| + |x_{t_i}(x_i) - x_t(x)|$$
$$\leq 2(F_n + F)\sigma_n^\beta + |z_{t_i}^n(x_i) - x_{t_i}(x_i)|.$$

It follows that

$$\sup_{0 \leq t \leq 1} \sup_{|x| \leq R} |z_t^n(x) - x_t(x)|^p \leq C_p \bigg\{ (F_n^p + F^p)\sigma_n^{\beta p} + \sup_{1 \leq i \leq N_n} |z_{t_i}^n(x_i) - x_{t_i}(x_i)|^p \bigg\}$$
$$\leq C_p \bigg\{ (F_n^p + F^p)\sigma_n^{\beta p} + \sum_{1 \leq i \leq N_n} |z_{t_i}^n(x_i) - x_{t_i}(x_i)|^p \bigg\}$$

with a constant $C_p$ depending only on $p$. Therefore, for another such constant $\hat{C}_p > 0$, we have

$$\mathbf{E}\bigg( \sup_{0 \leq t \leq 1} \sup_{|x| \leq R} |z_t^n(x) - x_t(x)|^p \bigg) \leq \hat{C}_p \sigma_n^{\beta p} + N_n \varepsilon_n$$
$$\leq \hat{C}_p \sigma_n^{\beta p} + C\bigg(\frac{1}{\sigma_n}\bigg)^{d+1} \cdot \varepsilon_n.$$

Now taking $\sigma_n = \varepsilon_n^{1/2(d+1)}$ gives the result (3.19). □

Due to hypothesis (H), Theorem 3.4 finally implies Theorem A, following a procedure in Chapter V of [4].

S. FANG
INSTITUT DE MATHÉMATIQUES
UNIVERSITÉ DE BOURGOGNE
B.P. 47870
21078 DIJON
FRANCE
E-MAIL: fang@u-bourgogne.fr

P. IMKELLER
INSTITUT FÜR MATHEMATIK
HUMBOLDT-UNIVERSITÄT ZU BERLIN
UNTER DEN LINDEN 6
10099 BERLIN
GERMANY
E-MAIL: imkeller@mathematik.hu-berlin.de

T. ZHANG
SCHOOL OF MATHEMATICS
UNIVERSITY OF MANCHESTER
OXFORD ROAD
MANCHESTER M13 9PL
ENGLAND
E-MAIL: tzhang@maths.man.ac.uk